\documentclass{amsart}
\usepackage{amsmath,amsfonts,amsthm}
\usepackage{amssymb}

\theoremstyle{plain}
\newtheorem{thm}{Theorem}[section]
\newtheorem{lem}[thm]{Lemma}

\theoremstyle{definition}
\newtheorem{defi}[thm]{Definition}

\theoremstyle{remark}
\newtheorem{rmk}[thm]{Remark}

\numberwithin{equation}{section}

\begin{document}

\title{On endpoint regularity criterion of the 3D Navier-Stokes equations}

\author{Zhouyu Li}
\address[Zhouyu Li]{School of Mathematics and Statistics, Northwestern Polytechnical University, Xi'an, Shaanxi 710129, P. R. China}
\email{zylimath@163.com}

\author{Daoguo Zhou}
\address[Daoguo Zhou]{School of Mathematics and Information Science, Henan Polytechnic University, Jiaozuo, Henan 454000, P. R. China}
\email{daoguozhou@hpu.edu.cn}


\date{\today}

\begin{abstract}
Let $(u, \pi)$ with $u=(u_1,u_2,u_3)$ be a suitable weak solution of the three dimensional Navier-Stokes equations in $\mathbb{R}^3\times (0, T)$. Denote by $\dot{\mathcal{B}}^{-1}_{\infty,\infty}$
the closure of $C_0^\infty$ in $\dot{B}^{-1}_{\infty,\infty}$. We prove that if
 $u\in L^\infty(0, T; \dot{B}^{-1}_{\infty,\infty})$,
$u(x, T)\in \dot{\mathcal{B}}^{-1}_{\infty,\infty})$, and
$u_3\in L^\infty(0, T; L^{3, \infty})$
or $u_3\in L^\infty(0, T; \dot{B}^{-1+3/p}_{p, q})$ with $3<p, q< \infty$,
then $u$ is smooth in $\mathbb{R}^3\times (0, T]$. 
Our result improves a previous result established by Wang and Zhang  [Sci. China Math. 60, 637-650 (2017)].
\end{abstract}

\maketitle

\section{Introduction}

In this paper, we study the incompressible Navier-Stokes equations in $\mathbb{R}^3\times (0, T)$
\begin{equation}\label{N-S}
    \begin{cases}
       \partial_{t}u-\Delta u +u\cdot \nabla u +\nabla \pi=0,\\
       \text{div}\;u=0,\\
       u|_{t=0}=u_0(x),
    \end{cases}
    \end{equation}
where $u(x, t)=(u_1,u_2, u_3)$ denotes the velocity of the fluid, and $\pi(x, t)$ represents
the pressure.

In the pioneering works \cite{Hopf1951, Leray1934}, Leray and Hopf proved the global existence
of weak solutions with finite energy. However, the uniqueness and regularity of weak solutions still remains open. There exists many conditional regularity results of the three dimensional Navier-Stokes
equations. The most well-known one is due to  Serrin \cite{Serrin1962} (see also Struwe \cite{Struwe1988}), which states that if the weak solution $u$ satisfies
\[
      u\in L^q(0, T; L^p(\mathbb{R}^3)) \quad \mbox{with} \quad \frac{2}{q}+\frac{3}{p}\leq 1, \quad 3<p\leq\infty,
      \]
then $u$ is regular. The limiting case $u\in{L^\infty(0, T; L^3(\mathbb{R}^3))}$ is solved by
Escauriaza, Seregin and \v{S}ver\'ak \cite{ESS-2003}  using blow-up analysis and backward
uniqueness for heat equations.

In term of the chain of critical spaces
\[L^3(\mathbb{R}^3)\hookrightarrow L^{3,q}(\mathbb{R}^3)\hookrightarrow \dot{B}^{-1+3/p}_{p, q}(\mathbb{R}^3)(3<p,q<\infty)\hookrightarrow \mathrm{BMO^{-1}}(\mathbb{R}^3)\hookrightarrow\dot{B}_{\infty, \infty}^{-1}(\mathbb{R}^3),\]
it is natural to extend the regularity condition $u\in{L^\infty(0, T; L^3(\mathbb{R}^3))}$ to $u\in{L^\infty(0, T; X )}$ with $X$ being one of the above spaces. Here,  a Banach space $X$ is called critical if  we have $\|\lambda u(\lambda x)\|_X=\|u\|_X$ for $u\in X$. Many works have been devoted to this research direction.
In particular, Gallagher et al.  \cite{Gallagher2016} proved the case
$u\in L^\infty(0, T; \dot{B}^{-1+3/p}_{p, q}(\mathbb{R}^3))$ with
$3<p, q<\infty$.
Recently, Wang and Zhang \cite{WangZhang2017} further generalized the result of  Gallagher, Koch and Planchon to the condition:
\[
u_3\in L^\infty(0, T; \dot{B}^{-1+3/p}_{p, q}(\mathbb{R}^3)) \quad \mbox{for} \quad 3<p, q<\infty,
\]
and
\[u\in L^\infty(0, T;\mathrm{BMO^{-1}} (\mathbb{R}^3)) \quad\mbox{with} \quad u(T)\in \mathrm{VMO^{-1}(\mathbb{R}^3)}.\] For more related results, see \cite{Alb2018,AlbBarker2019,BarkerSeregin2017,DongDu2009,LiWang2019, Phuc2015,Seregin2012,Tao2019} and references therein.

In this paper, we aim at improving  Wang and Zhang's results to the case
$u\in L^\infty(0, T; \dot{B}^{-1}_{\infty,\infty}(\mathbb{R}^3))$. Our main results read as follows.

\begin{thm}\label{Main-thm}
Let $(u, \pi)$ be a suitable weak solution of \eqref{N-S} in $\mathbb{R}^3\times (0, T)$. If
 $u\in L^\infty(0, T; \dot{B}^{-1}_{\infty,\infty}(\mathbb{R}^3))$,
$u(x, T)\in \dot{\mathcal{B}}^{-1}_{\infty,\infty}(\mathbb{R}^3)$, and $u_3$ satisfies
\[u_3\in L^\infty(0, T; L^{3, \infty}(\mathbb{R}^3)),\]
or
\[u_3\in L^\infty(0, T; \dot{B}^{-1+3/p}_{p, q}(\mathbb{R}^3))\quad  \text{with} \quad 3<p, q< \infty,\]
then $u$ is smooth in $\mathbb{R}^3\times (0, T]$. Here,  $\dot{\mathcal{B}}^{-1}_{\infty,\infty}(\mathbb{R}^3)$
denotes the closure of $C_0^\infty(\mathbb{R}^3)$ in $\dot{B}^{-1}_{\infty,\infty}(\mathbb{R}^3)$.
\end{thm}
\begin{rmk}
Considering the inclusion relationship that for $3<p,q<\infty$,
  \[
    L^{3}(\mathbb{R}^3)\hookrightarrow L^{3,q}(\mathbb{R}^3)\hookrightarrow\dot{B}^{-1+3/p}_{p, q}(\mathbb{R}^3)\hookrightarrow VMO^{-1}(\mathbb{R}^3)\hookrightarrow \mathcal{\dot{B}}^{-1}_{\infty,\infty}(\mathbb{R}^3),\]
Theorem \ref{Main-thm} improves previous results due to Escauriaza, Seregin and \v{S}ver\'ak \cite{ESS-2003}, Phuc \cite{Phuc2015}, Gallagher, Koch and Planchon  \cite{Gallagher2016} and
  Wang and Zhang \cite{WangZhang2017}.
\end{rmk}

\begin{rmk}
 Denote by $\dot{\mathbb{B}}^{-1}_{\infty,\infty}(\mathbb{R}^3)$ the closure of $L^\infty(\mathbb{R}^3)$ in $\dot{B}^{-1}_{\infty,\infty}(\mathbb{R}^3)$. It holds that  $\dot{B}^{-1}_{\infty,q}(\mathbb{R}^3)\hookrightarrow\dot{\mathbb{B}}^{-1}_{\infty,\infty}(\mathbb{R}^3)$ for $1\leq q< \infty$. Tobias Barker pointed out that $\dot{\mathcal{B}}^{-1}_{\infty,\infty}(\mathbb{R}^3)$ in Theorem \ref{Main-thm} can be replaced with $\dot{\mathbb{B}}^{-1}_{\infty,\infty}(\mathbb{R}^3)$. Namely, we have for every $f\in \dot{\mathbb{B}}^{-1}_{\infty,\infty}(\mathbb{R}^3)$,
\begin{equation}\label{vanish}
  \lim_{\lambda\to 0}\lambda f(\lambda \cdot)=0 \quad \text{in}\quad \mathcal{D}'.
\end{equation}
\end{rmk}

Our proof of Theorem \ref{Main-thm} is based on the scheme developed by Escauriaza, Seregin and \v{S}ver\'ak \cite{ESS-2003}, which consists of blow up analysis and backward uniqueness of heat equations. We need to show two properties of the limit of scaled solutions: it vanishes at the last moment and has some spatial decay.
We establish the first property by combining the facts that $\dot{B}^{-1}_{\infty,\infty}(\mathbb{R}^3)$ is invariant under the Naiver-Stokes scaling and $C_0^\infty(\mathbb{R}^3)$ is dense in $\dot{\mathcal{B}}^{-1}_{\infty,\infty}(\mathbb{R}^3)$. To get spatial decay for the limit function,
we take advantage of the property of involved function spaces --
$L^p$ is dense in $\dot{B}^{-1+3/p}_{p, q}(\mathbb{R}^3)(3<p,q<\infty)$, and $L^{3, \infty}(\mathbb{R}^3)\subset L^2(\mathbb{R}^3)+L^4(\mathbb{R}^3)$. In contrast, to show decay of functions in $\dot{B}^{-1+3/p}_{p, q}(\mathbb{R}^3)(3<p,q<\infty)$, Wang and Zhang \cite{WangZhang2017} used the fact that $C_0^\infty$ is dense in $\dot{B}^{-1+3/p}_{p, q}(\mathbb{R}^3)(3<p,q<\infty)$, which is nontrivial (see Theorem 3.5 in Chapter 3 of Triebel \cite{Triebel}).

The rest of this paper is organized as follows. In Section 2, we recall the definitions of weak solutions and suitable weak solutions, as well as some useful
lemmas. Section 3 is devoted to the proof of main result.

\section{Preliminaries}

In this section we state definitions of Leray-Hopf weak solutions and suitable weak solutions,
introduce some notations, and collect some useful lemmas.

First we recall the definitions of Leray-Hopf weak solutions \cite{Hopf1951,Leray1934}  and suitable weak solutions  \cite{CKN1982} to the
Navier-Stokes equations.

\begin{defi}
A vector field $u$ is called a Leray-Hopf weak solution of
	\eqref{N-S} in $\mathbb{R} ^3\times (0, T)$ if
\begin{itemize}
  \item [(1)]
 $u\in L^\infty(0, T; L^2(\mathbb{R} ^3))\cap L^2(0, T; H^1(\mathbb{R} ^3))$;

  \item  [(2)]
  $u$ satisfies \eqref{N-S} in $\mathbb{R} ^3\times (0, T)$ in the weak sense that
	for any $\psi\in C_c^\infty(\mathbb{R}^3\times(0,T))$ such that $\text{div}\; \psi=0$,
  \[\int_0^T\int_{\mathbb{R} ^3}(-u\partial_t\psi+\nabla u\nabla \psi-uu\nabla\psi)dxdt=0.\]

  \item  [(3)]
  $u$ satisfies the energy inequality
  \begin{equation*}
	\int_{\mathbb{R} ^3}|u(x, t)|^2 dx +2\int_{0}^t\int_{\mathbb{R} ^3}|\nabla u|^2 dxds\\
 \leq \int_{\mathbb{R} ^3}|u_0(x)|^2 dx
	\end{equation*}
  for all a.e. $t\in [0, T]$.
\end{itemize}
\end{defi}

\begin{defi}
Let $\Omega$ be an open set in $\mathbb{R}^3$ and $T>0$. A pair $(u, \pi)$ is called a suitable weak solution of
\eqref{N-S} in $\Omega \times (0, T)$ if
\begin{itemize}
  \item [(1)]
$u\in L^\infty(0,T; L^2(\Omega))\cap L^2(0, T; H^1(\Omega))$,  $\pi\in L^{\frac{3}{2}}(\Omega \times (0, T))$;

\item [(2)]
 $(u, \pi)$ satisfies \eqref{N-S} in $\Omega \times (0, T)$ in the sense of distribution;

 \item [(3)]
 $(u, \pi)$ satisfies the local energy inequality
\begin{equation*}\label{2-2}
    \begin{split}
       &\int_\Omega|u(x, t)|^2\varphi \, dx +2\int_0^t\int_\Omega|\nabla u|^2 \varphi \, dxds\\
       &\quad \leq \int_0^t\int_\Omega|u|^2(\partial_s \varphi+ \Delta\varphi)+
       u\cdot\nabla \varphi (|u|^2+2\pi) \, dxds
    \end{split}
    \end{equation*}
for all a.e. $t\in [0, T]$ and all  $\varphi\in C_c^\infty(\mathbb{R}^3\times\mathbb{R})$ such that $\varphi\geq 0$
in $\Omega \times (0, T)$.
\end{itemize}
\end{defi}

\begin{rmk}
If $u$ is a Leray-Hopf weak solution and $u\in L^\infty(0, T; \dot{B}^{-1}_{\infty,\infty}(\mathbb{R}^3))$,
by Lemma \ref{Lem-int}, then
$u\in L^4(\mathbb{R}^3\times (0, T))$, which means that $u$ satisfies the local energy inequality and is a suitable
weak solution.
\end{rmk}

We now fix some notations. Let $(u, \pi)$ be a solution to the Navier-Stokes equations \eqref{N-S}.
Then for $\lambda \in \mathbb{R}$ $(u_\lambda(x, t),\pi_\lambda(x, t))$
\[
u_\lambda(x, t)=\lambda u(\lambda x, \lambda^2 t),\quad\quad \pi_\lambda(x, t)=\lambda^2 \pi(\lambda x, \lambda^2 t),\]
also solves the Navier-Stokes equations \eqref{N-S}. For $z_0=(x_0, t_0)$, the following quantities are invariant
under the above scaling:
\[
A(u; r, z_0):= \sup_{-r^2+t_0\leq t\leq t_0}r^{-1}\int_{B_r(x_0)}|u(y, t)|^2dy, \;
C(u;r, z_0):= r^{-2}\int_{Q_r(z_0)}|u(y, s)|^3dyds,
\]
\[
E(u;r, z_0):= r^{-1}\int_{Q_r(x_0)}|\nabla u(y, s)|^2dyds,  \;
D(u;r, z_0):= r^{-2}\int_{Q_r(x_0)}|\pi(y, s)|^\frac{3}{2}dyds.
\]
We denote
\[
B(x_0, r):=\{x\in \mathbb{R}^3: |x-x_0|< R\},\quad B_r:=B(0, r);
\]
\[
Q(z_0, r):=B(x_0, r)\times (-r^2+t_0, t_0), \quad Q_r:=Q(0, r);
\]
\[
A(u;r):=A(u;r, 0),
\]
and so on.

We now recall the definitions of Littlewood-Paley decomposition and Besov space, and imbedding theorem  in Besov spaces. Then we collect some useful properties of Besov spaces.

\begin{defi}[Littlewood-Paley decomposition]
  Let $\phi$ be a smooth function with values in $[0, 1]$ such that
  $\phi \in C_0^\infty (\mathbb{R}^n)$ in the annulus 
  $\mathcal{C}:=\{\xi\in \mathbb{R}^n\colon \frac{3}{4}\leq|\xi|\leq \frac{8}{3}  \}$ satisfying
  \[
   \sum_{j\in \mathbb{Z}}\phi (2^{-j}\xi)=1, \quad \, \forall \xi\in \mathbb{R}^n\setminus\{0\}.
   \]
    For every $u\in S'$, we define the homogeneous dyadic blocks $\dot{\Delta}_{j}$ and the homogeneous low-frequency cutoff
    operators $S_j$ for all $j\in \mathbb{Z}$ as
\[
  \dot{\Delta}_j f =\phi(2^{-j}D) f = \mathcal{F}^{-1}(\phi(2^{-j}\xi)\hat{f}(\xi)).
 \]
\end{defi}

\begin{defi}
    Let $P$ be the set of polynomials. The Besov space
    $\dot B^{s}_{p,q}(\mathbb{R}^3)$ with  $s\in \mathbb{R}$, $1\leq p,q\leq \infty$, consists of $f\in \mathcal{S}' (\mathbb{R}^3)/P$ satisfying
    \[\displaystyle{\|f\|_{\dot B^s_{p,q}}:= \big\| 2^{js}\| \dot{\Delta}  f \|_{L^{p}(\mathbb{R}^3)} \big\|_{\ell^q} <\infty}.\]
\end{defi}
  \begin{lem}\label{besovimb}
    Let $1\leq p_1\leq p_2\leq \infty$, and $1\leq r_1\leq r_2\leq \infty$, $s\in \mathbb{R}$. Then, we have
    \[\dot{B}^{s}_{p_1,r_1}\hookrightarrow \dot{B}^{s+d/{p_2}-d/{p_1}}_{p_2,r_2}.\]
  \end{lem}
\begin{rmk}
  Lemma \ref{besovimb} can be found in Chapter 2 in Bahouri et.al. \cite{BahChemin2011}.
\end{rmk}

The following improved Galirado-Nirenberg inequality (Ledoux \cite{Ledoux2003}) and its local version (Seregin and the second author \cite{SereginZhou2020}) will be used in the proof of main theorem.
\begin{lem}\label{Lem-int}
 (1) If $f\in \dot{B}^{-1}_{\infty,\infty}(\mathbb{R}^3)\cap H^1(\mathbb{R}^3)$, then we have
    \[
      \|f\|_{L_4(\mathbb{R}^3)}\leq c\|f\|^{\frac{1}{2}}_{\dot{B}^{-1}_{\infty,\infty}(\mathbb{R}^3)}\|\nabla f\|^{\frac{1}{2}}_{L^2(\mathbb{R}^3)}.
    \]
(2) If $f\in \dot{B}^{-1}_{\infty,\infty}(\mathbb{R}^3)\cap H^1(B_{2r}(x_0))$, then we have
    \[
      \|f\|_{L_4(B_r(x_0))}\leq c\|f\|^{\frac{1}{2}}_{\dot{B}^{-1}_{\infty,\infty}(\mathbb{R}^3)}\Big(\frac{1}{r}\|f\|_{L^2(B_{2r}(x_0))}+\|\nabla f\|_{L^2(B_{2r}(x_0))}\Big)^{\frac{1}{2}}.
      \]
\end{lem}

We recall the definition of weak Lebesgue space, as well as decomposition of functions in weak Lesbegue spaces (see \cite{BarkerSeregin2017}).
\begin{defi}
  The weak Lesbegue space $L^{3,\infty}(\mathbb{R}^3)$ consists of local integrable functions $f$ satisfying
  \[\|f\|_{L^{3,\infty}(\mathbb{R}^3)}=\sup_{\lambda>0}\lambda|\{x\in \mathbb{R}^3\colon |u|>\lambda\}|^{\frac 1 3}<\infty.\]
  \end{defi}
\begin{lem}\label{Lem-4}
Let $1<t<r<s\leq\infty$, and  $f\in L^{r, \infty}(\mathbb{R}^3)$.
Then we have $f=f_1+f_2$ for some $f_1\in L^s(\mathbb{R}^3)$ and $f_2\in L^t(\mathbb{R}^3)$, which satisfy
\begin{equation*}\label{1-1}
\|f_1\|_{L^s(\mathbb{R}^3)}\leq c(s, r)\|f\|_{L^{r, \infty}(\mathbb{R}^3)},
\end{equation*}
and
\begin{equation*}\label{1-2}
\|f_2\|_{L^t(\mathbb{R}^3)}\leq c(r, t)\|f\|_{L^{r, \infty}(\mathbb{R}^3)}.
\end{equation*}
\end{lem}

We also need a bound for scaled energy of Naiver-Stokes equations by Seregin and the second author \cite{SereginZhou2020}.
\begin{lem}\label{Lem-1}
Let $(u, \pi)$ be a suitable weak solution of \eqref{N-S} in $\mathbb{R}^3\times (0, T)$.
Moreover, it is supposed that
\[u\in L^\infty(0, T; \dot{B}^{-1}_{\infty, \infty}(\mathbb{R}^3)).\]
Then, for any $z_0\in \mathbb{R}^3\times(0, T]$, we have estimate
\begin{equation*}
  \begin{gathered}
\sup_{0<r<r_0}\big\{A(u;z_0, r)+ C(u;z_0, r)+ D(u;z_0, r)+ E(u;z_0, r)\big\}\\
\leq c\big[r_0^{\frac{1}{2}}+\|u\|^2_{L^\infty(0, T; \dot{B}^{-1}_{\infty, \infty}(\mathbb{R}^3))}+\|u\|^6_{L^\infty(0, T; \dot{B}^{-1}_{\infty, \infty}(\mathbb{R}^3))}\big],
 \end{gathered}
\end{equation*}
where $r_0\leq \frac{1}{2}\min \{1,t_0\}$ and $c$ depends on $C(z_0, 1)$ and $D(z_0, 1)$ only.
\end{lem}

We conclude this section by recalling the small energy regularity results due to Wang and Zhang \cite{WangZhang2017}.
\begin{lem}\label{Lem-2}
Let $(u, \pi)$ be a suitable weak solution of \eqref{N-S} in $Q_1$. If $u$ satisfies
\[\sup_{0<r<1}\big\{A(u;r)+E(u;r)\big\}\leq M,\]
where $M>0$, then there exists a positive
constant $\varepsilon$ depending on $M$ such that if
\[\frac{1}{r_*^2}\int_{Q_{r_*}}|u_3|^3dxdt\leq \varepsilon,\]
for some $r_*$ with $0<r_*< \min \{\frac{1}{2}, (C(u;1)+D(\pi;1))^{-2}\}$, then $(0, 0)$ is a regular point.
\end{lem}

\begin{lem}\label{Lem-3}
Let $(u, \pi)$ be a suitable weak solution of \eqref{N-S} in $Q_r$. If $(u, \pi)$ satisfies
\[\frac{1}{r^2}\int_{Q_r}|u|^3+|\pi|^{3/2}dxdt\leq M,\]
where  $M>0$, then there exists a positive constant $\varepsilon$ depending on $M$ such that if
\[\frac{1}{r^2}\int_{Q_r}|u_3|^3dxdt\leq \varepsilon,\]
then $(0, 0)$ is a regular point.
\end{lem}

\section{Proof of Theorem \ref{Main-thm}}

Since the Naiver-Stokes equations are translation and scaling invariant, thus it is sufficient to prove our main results in the domain
$\mathbb{R}^3\times [-1, 0)$. The proof is based on blow-up analysis and backward uniqueness of parabolic equations developed in Escauriaza, Seregin and \v{S}ver\'ak \cite{ESS-2003}. We argue by contradiction.  Without loss of generality, we assume
that $(0, 0)$ is a singular point of $u$.

\begin{proof}
{\bf Step 1: Blow-up analysis}

By the assumption, we have
\[\|u\|_{L^\infty(-1, 0; \dot{B}^{-1}_{\infty, \infty}(\mathbb{R}^3))}\leq c.\]
From Lemma \ref{Lem-1}, we get, for $z_0\in B_{1/2}\times (-1/4,0)$ and $0< r < \frac{1}{2}$,
\begin{equation}\label{boundeu}
A(u;z_0,r)+C(u;z_0,r)+D(u;z_0,r)+E(u;z_0,r) \leq c(C(u;1), D(u;1)).
\end{equation}

Since $(0, 0)$ is a singular point, Lemma \ref{Lem-2} ensures that there exists a sequence $R_k$
such that $R_k\rightarrow 0$ as $k\rightarrow +\infty$ and
\begin{equation}\label{3-2}
R_k^{-2}\int_{Q_{R_k}}|u_3(x, t)|^3 \, dxdt \geq \varepsilon.
\end{equation}

Define
\[
u^k(y, s)=R_k u(R_ky, R_k^2s), \quad \pi^k(y, s)=R_k^2 \pi(R_ky, R_k^2s),
\]
where $(y, s)\in \mathbb{R}^3\times (-\frac{1}{R_k^2}, 0)$. Then the pair $(u^k, \pi^k)$ is still a suitable weak solution to ~\eqref{N-S}.

Since $A(u; r)$, $C(u; r)$, $D(\pi; r)$ and $E(u; r)$ are invariant under the Navier-Stoeks scaling, we obtain that, for any $a>0$ and $z_0=(x_0,t_0)\in \mathbb{R}^3\times (-\infty, 0]$,
\begin{equation*}
\begin{aligned}
&A(u^k;z_0,a)+C(u^k;z_0,a)+D(\pi^k;z_0,a)+E(u^k; z_0,a) \\
&=A(u^k;z_0^k,r_ka)+C(u^k;z_0^k,r_ka)+D(\pi^k;z_0^k,r_ka)+E(u^k; z_0^k,r_ka),
\end{aligned}
\end{equation*}
where $z_0^k=(r_kx_0,r_k^2t_0)$. From \eqref{boundeu}, we get that for enough large $k$,
\begin{equation*}
A(u^k;z_0,a)+C(u^k;z_0,a)+D(\pi^k;z_0,a)+E(u^k; z_0,a)\leq c.
\end{equation*}

We also have
\begin{equation*}\label{3-11}
 \int_{Q_1} |u_{3}^k(x, t)|^3 \, dxdt= R_k^{-2}\int_{Q_{R_k}}|u_3(x, t)|^3\, dxdt \geq \varepsilon
\end{equation*}
for all $k\in\mathbb{N}$.

By interpolation between $A(u^k; a)$ and $E(u^k; a)$, we get $u^k\in L^6_tL_x^{18/7}(Q_a)$. Then  we have by H\"older's inequality that $u^k\cdot \nabla u^k\in L_t^\frac{3}{2}L_x^\frac{9}{8}(Q_a)$. Appealing to the linear
Stokes estimate, we deduce that
\[|\partial_t u^k|+|\Delta u^k|+ |\nabla \pi^k|\in L_t^\frac{3}{2}L_x^\frac{9}{8}(Q_a).\]

Applying the  Aubin-Lions lemma, we can extract a subsequence, still denoted by
$(u^k,\pi^k)$, such that $(u^k,\pi^k)$  converges weakly to some limit functions
$(v, \pi')$,  for any $a > 0$,
\begin{equation}\label{covlim}
\begin{gathered}
u^k \rightharpoonup v \quad \text{ in} \quad L^\infty(-a^2,0; \dot{B}^{-1}_{\infty, \infty}(\mathbb{R}^3)),\\
u^k \rightharpoonup v \quad \text{ in} \quad L^\infty(-a^2, 0; L^2(B_a)),\\
u^k \rightarrow v \quad \text{strongly in}  \quad C([-a^2, 0]; L^\frac{9}{8} (B_a)),\\
\nabla u^k \rightharpoonup \nabla v \quad \text{in} \quad L^2(Q_a),\\
\pi^k \rightharpoonup \pi'\quad \mbox{in} \quad L^\frac{3}{2}(Q_a).
\end{gathered}
\end{equation}
Since $u^k$ is uniformly bounded in $L^{10/3}(Q_a)$ by interpolation between $A(u^k; a)$ and $E(u^k; a)$, we get by H\"older's inequality,
\begin{equation}\label{covL3}
u^k \rightarrow v \quad \text{strongly in} \quad L^3(Q_a).
\end{equation}
Furthermore, we have in case  $u_3\in L^\infty(-1, 0; L^{3,\infty}(\mathbb{R}^3))$ that
\begin{equation*}
u_3^k \rightharpoonup v_3 \quad  \text{ in}  \quad L^\infty(-a^2, 0; L^{3,\infty}(\mathbb{R}^3)),
\end{equation*}
or in case $u_3\in L^\infty(-1, 0; \dot{B}^{-1+3/p}_{p, q}(\mathbb{R}^3))$
that
\begin{equation*}
u_3^k \rightharpoonup v_3\quad  \text{ in}  \quad L^\infty(-a^2, 0; \dot{B}^{-1+3/p}_{p, q}(\mathbb{R}^3)).
\end{equation*}

The above convergence implies that $(v,\pi')$ satisfy the
Navier-Stokes equations in $\mathbb{R}^3\times (-\infty,0)$.
Moreover, due to lower semi-continuity of norm, it holds that for any $z_0\in \mathbb{R}^3\times (-\infty,0)$ and $a>0$,
\begin{equation}\label{bounde}
A(v; z_0,a)+E(v; z_0,a)+C(v; z_0,a)+D(\pi'; z_0,a) \leq c,
\end{equation}
\begin{equation}
\|v_3\|_{L^\infty(-a^2, 0; \dot{B}^{-1+3/p}_{p, q}(\mathbb{R}^3))}\leq c, \text{ if } u_3\in L^\infty(-a^2, 0; \dot{B}^{-1+3/p}_{p, q}(\mathbb{R}^3));
\end{equation}
\begin{equation}
\|v_3\|_{L^\infty(-a^2, 0; L^{3, \infty}(\mathbb{R}^3))}\leq c, \text{ if } u_3\in L^\infty(-a^2, 0; L^{3,\infty}(\mathbb{R}^3)).
\end{equation}
Thanks to \eqref{covL3},  we find
\begin{equation}\label{contra}
\int_{Q_1} |v_{3}(x, t)|^3 \, dxdt \geq \varepsilon.
\end{equation}

{\bf Step 2:  Prove that the limit function $v(x, 0)=0$ in $\mathbb{R}^3$}

Since $u(x, 0)\in \mathcal{\dot{B}}^{-1}_{\infty, \infty}(\mathbb{R}^3)$,
where $\mathcal{\dot{B}}^{-1}_{\infty, \infty}(\mathbb{R}^3)$ is the closure of $C_0^\infty(\mathbb{R}^3)$ in $\dot{B}^{-1}_{\infty,\infty}(\mathbb{R}^3)$, thus for any $\varepsilon > 0$, there exists a function
$U^*(x, 0) \in C_0^\infty (\mathbb{R}^3)$ such that
$\|u(x,0)-U^*\|_{\dot{B}^{-1}_{\infty, \infty}(\mathbb{R}^3)} < \varepsilon$. Then for any $\varphi \in C_c^\infty (B_a)$ with $a > 0$, we have
\begin{equation*}\label{4-1}
  \begin{split}
|\int_{B_a}v(x, 0)\varphi(x) \, dx|&\leq |\int_{B_a}(v(x, 0)- u^k(x, 0))\varphi(x) \, dx| + |\int_{B_a}u^k(x, 0)\varphi(x) \, dx|\\
&=I_1+I_2.
 \end{split}
\end{equation*}

For $I_1$, the convergence in \eqref{covlim} shows that
\begin{equation*}\label{4-2}
I_1\leq \int_{B_a}|v(x, 0)- u^k(x, 0)| \, dx\, \rightarrow 0 \quad \mbox{as} \quad k \rightarrow+\infty .
\end{equation*}

For $I_2$, we obtain
\begin{equation*}
  \begin{split}
I_2&= |\int_{B_a}u^k(x, 0)\varphi(x) \, dx|\\
  &= R_k|\int_{B_a}u(R_kx,0)\varphi(x) \, dx|\\
  &\leq R_k|\int_{B_a} (u(R_kx,0)-U^*(R_kx))\varphi(x) \, dx|+ R_k|\int_{B_a} U^*(R_kx)\varphi(x) \, dx|\\
  &\leq cR_k\|u(R_kx,0)-U^*(R_kx)\|_{\dot{B}_{\infty, \infty}^{-1}(\mathbb{R}^3)}+ R_k\int_{B_a} |U^*(R_kx)\varphi(x)| \, dx\\
   &\leq c\|U(x)-U^*(x)\|_{\dot{B}_{\infty, \infty}^{-1}(\mathbb{R}^3)}+ R_k\int_{B_a} |U^*(R_kx)\varphi(x)| \, dx\\
  &\leq c\varepsilon,
 \end{split}
\end{equation*}
as $k \rightarrow+\infty$, where we used the fact that $U^*$ is continuous at $0$.

Thus we conclude  that $v(x, 0)=0$.

{\bf Step 3: Spatial  Decay}

(1){ \bf Case $u_3\in L^\infty(0, T; L^{3, \infty}(\mathbb{R}^3))$}

Using Lemma \ref{Lem-4}, we decompose $v_3(\cdot, t)\in L^{3, \infty}(\mathbb{R}^3)$  as
\begin{equation*}\label{6-2}
v_3(\cdot, t)=\bar{v}_3(\cdot, t)+\tilde{v}_3(\cdot, t),
\end{equation*}
where $\bar{v}_3(\cdot, t)\in L^2(\mathbb{R}^3)$,
$\tilde{v}_3(\cdot, t)\in L^4(\mathbb{R}^3)$,
and
\begin{equation*}\label{6-3}
\|\bar{v}_3(\cdot, t)\|_{L^2}\leq c\|v_3(\cdot, t)\|_{L^{3, \infty}(\mathbb{R}^3)},
\end{equation*}
\begin{equation*}\label{6-4}
\|\tilde{v}_3(\cdot, t)\|_{L^4}\leq c\|v_3(\cdot, t)\|_{L^{3, \infty}(\mathbb{R}^3)}.
\end{equation*}
Then we have $\bar{v}_3(x, t)\in L^\infty(-a^2, 0; L^2(\mathbb{R}^3))$ and
$\tilde{v}_3(x, t)\in L^\infty(-a^2, 0; L^4(\mathbb{R}^3))$.

Then it follows that
\begin{equation}\label{decayw}
\begin{split}
\int_{Q_1(z_0)}|v_3|^2 \, dxdt&\leq \int_{Q_1(z_0)}|\bar{v}_3|^2 \, dxdt + \int_{Q_1(z_0)}|\tilde{v}_3|^2 \, dxdt\\
&\rightarrow 0 \quad \text{as} \quad |z_0|\rightarrow \infty.
\end{split}
\end{equation}

(2) {\bf Case $u_3\in L^\infty(0, T; \dot{B}^{-1+3/p}_{p, q}(\mathbb{R}^3))$, $3< p, q<\infty$}

Let $v_3^N=\sum_{i=-N}^{N}\dot{\Delta}_i v_3$. We first show a functional property of Besov space: for any $\beta\in [1, \infty)$ and $T>0$, it holds that
\begin{equation}\label{approxbes}
  \|v_3-v^N_3\|_{L^\beta(-T, 0; \dot{B}^{-1+3/p}_{p, q}(\mathbb{R}^3))}\rightarrow 0 \quad \mbox{as} \quad N\rightarrow \infty.
\end{equation}
In fact, from the definition of Besov space, we have for any $t\in (-T,0]$,
\begin{equation*}\label{5-2a}
\lim_{N\rightarrow+\infty}\|v_3(t)-v^N_3(t)\|_{\dot{B}^{-1+3/p}_{p, q}(\mathbb{R}^3)}=0.
\end{equation*}
Since 
\begin{equation*}
\|v^N_3(t)\|_{\dot{B}^{-1+3/p}_{p, q}(\mathbb{R}^3)}\leq c \|v_3(t)\|_{\dot{B}^{-1+3/p}_{p, q}(\mathbb{R}^3)},
\end{equation*}
then \eqref{approxbes} is a consequence of Lebesgue's dominated convergence theorem.

By H\"older's inequality and Lemma \ref{Lem-int}, noting that
$\|v_3\|_{\dot{B}^{-1}_{\infty, \infty}(\mathbb{R}^3)}\leq c\|v_3\|_{\dot{B}^{-1+3/p}_{p, q}(\mathbb{R}^3)}$, we obtain
\begin{equation}\label{estiv3L2}
  \begin{aligned}
\int_{Q_1(z_0)}|v_3|^2 \, dxdt
&\leq \int_{Q_1(z_0)}|v_3-v_3^N|^2 \, dxdt+\int_{Q_1(z_0)}|v_3^N|^2 \, dxdt\\
&\leq  c\|v_3-v_3^N\|_{L^2(-T, 0; \dot{B}^{-1}_{\infty, \infty}(\mathbb{R}^3))}
\|v_3-v_3^N\|_{L^2_tH^1_x(Q_2(z_0))}+c\|v^N_3\|^2_{L^p(Q_1(z_0))}\\
&\leq c\|v_3-v_3^N\|_{L^2(-T, 0; \dot{B}^{-1+3/p}_{p, q}(\mathbb{R}^3)}
\big(\|v_3\|_{L^2_tH^1_x(Q_2(z_0))}+ c\|v_3^N\|_{L^2_tW^{1,p}_x(Q_2(z_0))}\big)\\
&\quad + c\|v^N_3\|^2_{L^p(Q_1(z_0))}.
\end{aligned}
\end{equation}

In \eqref{estiv3L2}, using \eqref{approxbes} and the fact that $\|v^N_3\|_{W^{1,p}(\mathbb{R}^3)}\leq c(N)\|v_3\|_{
  \dot{B}^{-1+3/p}_{p,q}(\mathbb{R}^3)}$, sending first $z_0$ to $\infty$ and then $N$ to $\infty$, we arrive at
\begin{equation}\label{decaybes}
\int_{Q_1(z_0)}|v_3| \, dxdt\rightarrow 0 \quad \mbox{as} \quad |z_0|\rightarrow \infty.
\end{equation}

{\bf Step 4: Backward uniqueness and unique continuation}

Spatial decay results in Step 3 and small regularity criterion in Lemma \ref{Lem-3} ensure that there exists a constant $R>0$ such that
\begin{equation}\label{regdecay}
|v(x, t)|+|\nabla v(x, t)|\leq c,
\end{equation}
for $(x, t)\in (\mathbb{R}^3\backslash B_R)\times (-T, 0)$.

Let $\omega=\nabla\times v$.  Since $v(x,0)=0$, we get $w(x,0)=0$. Moreover, as
$\omega$ satisfies
\begin{equation*}
\partial_t\omega-\Delta\omega=-v\cdot\nabla \omega+\omega\cdot\nabla v,
\end{equation*}
we find that by \eqref{regdecay}
\begin{equation*}
|\partial_t\omega-\Delta\omega|\leq c(|\omega|+|\nabla \omega|) \quad \text{in}\quad   (\mathbb{R}^3\backslash B_R)\times (-T, 0).
\end{equation*}
Applying backward uniqueness of parabolic operator \cite{ESS-2003} yields that
\begin{equation*}
\omega(x, t)=0,\quad (x, t)\in (\mathbb{R}^3\backslash B_R)\times (-T, 0).
\end{equation*}
By unique continuation argument as in \cite{ESS-2003}, we see that
\begin{equation}\label{vor0}
\omega(x, t)=0 \quad \text{in}\quad  \mathbb{R}^3\times (-T, 0).
\end{equation}

Combining the incompressible condition $\text{div}\; v=0$ and \eqref{vor0}, we find that $\Delta v(\cdot, t)=0$ in $\mathbb{R}^3$.
Then
we get  from the energy bound \eqref{bounde} that $v\in L^\infty(-T, 0; L^\infty(\mathbb{R}^3))$.

The Liouville Theorem for harmonic functions implies that $v(\cdot, t)$ is constant.  Due to the spacial decay bound in Step 3, it follows that $v_3(\cdot, t)=0$ for any $t\in(-T, 0)$,
which contradicts with \eqref{contra}. This completes  the proof of Theorem \ref{Main-thm}.
\end{proof}

\section*{Acknowledgments}
Z. Li was partially supported by the National Natural Science Foundation of China (No. 11601423) and the Natural Science Foundation of Shaanxi Province (No. 2020JQ-120). D. Zhou was partially supported by  the National Natural
Science Foundation of China (No. 11971446). D. Zhou wishes to thank Dallas Albritton, Tobais Barker and Wendong Wang for helpful discussions.


\begin{thebibliography}{99}

\bibitem{BahChemin2011}
H. Bahouri, J.Y. Chemin, R. Danchin, Fourier Analysis and Nonlinear Partial Differential Equations, Grundlehren der Mathematischen Wissenschaften, vol. 343, Springer, Heidelberg, 2011. xvi+523 pp.

\bibitem{Alb2018}
D. Albritton, Blow-up criteria for the Navier-Stokes equations in non-endpoint critical Besov spaces. {\it Anal. PDE,} {\bf11}  (2018), no. 6, 1415-1456.

\bibitem{AlbBarker2019}
D. Albritton, T. Barker, Global weak Besov solutions of the Navier-Stokes equations and applications.  {\it Arch. Ration. Mech. Anal.,} 232 (2019), no. 1, 197-263.

\bibitem{BarkerSeregin2017}
T. Barker, G. Seregin, A necessary condition of potential blowup for the Navier--Stokes
system in half-space. {\it Math. Ann.,} {\bf 369} (2017), no. 3-4, 1327--1352.

\bibitem{CKN1982}
L. Caffarelli, R. Kohn, L. Nirenberg, Partial regularity of suitable weak solutions of the Navier--Stokes equations.
{\it Comm. Pure Appl. Math.,} {\bf 35} (1982), no. 6, 771--831.

\bibitem{DongDu2009}
H.Dong, D. Du, The Navier-Stokes equations in the critical Lebesgue space. {\it Comm. Math. Phys.,} {\bf 292} (2009), no. 3, 811--827.

\bibitem{ESS-2003}
L. Escauriaza, G. Seregin, V. \v{S}ver\'ak, $L^{3, \infty}$ -- solutions of Navier--Stokes equations and backward uniqueness. (Russian) {\it Uspekhi Mat. Nauk,} {\bf 58} (2003), no. 2, 3--44.

\bibitem{Gallagher2016}
I. Gallagher, G. Koch, F. Planchon, Blow-up of critical Besov norms at a potential Navier--Stokes singularity.
{\it Comm. Math Phys.,} {\bf 343} (2016), no. 1, 39--82.

\bibitem{Hopf1951}
E. Hopf, \"Uber die Anfangswertaufgabe f\"ur die hydrodynamischen Grundgleichungen.
(German) {\it Math. Nachr.,} {\bf 4} (1951) 213--231.

\bibitem{Ledoux2003}
M. Ledoux, On improved Sobolev embedding theorems. {\it Math. Res. Lett.,} {\bf 10} (2003), no. 5-6,
659-669.

\bibitem{LiWang2019}
K. Li, B. Wang, Blowup criterion for Navier-Stokes equation in critical Besov space with spatial dimensions $d\geq 4$. {\it Ann. Inst. H. Poincar\'e Anal. Non Lin\'eaire,} {\bf 36} (2019), no. 6, 1679-1707.

\bibitem{Phuc2015}
N. Phuc, The Navier-Stokes Equations in Nonendpoint Borderline Lorentz Spaces. {\it J. Math. Fluid Mech.,} {\bf 17} (2015), no. 4, 741-760.

\bibitem{Leray1934}
J. Leray, Sur le mouvement d'un liquide visqueux emplissant l'espace.
{\it Acta Math.,} {\bf 63} (1934), no. 1, 193--248.

\bibitem{Serrin1962}
J. Serrin, On the interior regularity of weak solutions of the Navier-Stokes
equations. {\it Arch. Rational Mech. Anal.,} {\bf 9} (1962) 187--195.

\bibitem{Seregin2012}
G. Seregin,  A Certain Necessary Condition of Potential Blow up for Navier-Stokes Equations. {\it Comm. Math. Phys.,} {\bf 312} (2012), no. 3, 833-845.

\bibitem{SereginZhou2020}
G. Seregin, D. Zhou, Regularity of Solutions to the Navier-Stokes equations in $\dot{B}^{-1}_{\infty,\infty}.$
{\it J. Math. Sci.,}  {\bf 244} (2020), no. 6, 1003-1009.

\bibitem{Struwe1988}
M. Struwe, On partial regularity results for the Navier--Stokes equations.
{\it Comm. Pure Appl. Math.,} {\bf 41} (1988), no. 4, 437--458.

\bibitem{Tao2019}
T. Tao, Quantitative bounds for critically bounded solutions to the Navier-Stokes equations. {\it arXiv:1908.04958,} 2019.

\bibitem{Triebel}
H. Triebel, Tempered homogeneous function spaces. EMS Series of Lectures in Mathematics. European Mathematical Society (EMS), Z\"urich, 2015. xii+130 pp.

\bibitem{WangZhang2017}
W. Wang, Z. Zhang, Blow--up of critical norms for the 3-D Navier--Stokes equations. {\it Sci. China Math.,}
{\bf 60} (2017), no. 4, 637--650.

\end{thebibliography}
 \end{document}